\documentclass[11pt]{article}
\usepackage{amsmath,amssymb,amsthm,graphicx,psfrag,xspace}
\usepackage[bookmarks=true]{hyperref}

\def\rcs $#1: #2 ${\expandafter\def\csname rcs#1\endcsname {#2}}
\rcs $Date: 2008/06/17 19:09:31 $
\rcs $Revision: 1.27 $

\newtheorem{thm}{Theorem}[section]

\newtheorem{lem}[thm]{Lemma}
\newtheorem{prop}[thm]{Proposition}

\newcommand{\extraline}{\vskip\baselineskip}
\newcommand{\remark}{\noindent\textbf{Remark.}\xspace}

\renewcommand{\th}{^\text{th}}
\renewcommand{\Re}{\operatorname{Re}}
\renewcommand{\Im}{\operatorname{Im}}
\newcommand{\poly}{\mathcal{P}}
\newcommand{\harmconj}{\mathcal{H}}
\newcommand{\C}{\mathbb{C}}
\newcommand{\R}{\mathbb{R}}
\newcommand{\vk}{\vec v_k}
\newcommand{\Ik}{\vec I_k}
\newcommand{\Vk}{V_k}
\newcommand{\identity}{\mathrm{Id}}
\newcommand{\Z}{{\mathbb Z}}

\newcommand{\CdV}{MR1652692}
\newcommand{\IM}{MR1617177}
\newcommand{\CIM}{MR1657214}
\newcommand{\KW}{KW-polynomial}
\newcommand{\GM}{MR2150803}
\newcommand{\Gam}{MR1830078}
\newcommand{\Dub}{MR2253875}

\newcommand{\pu}[1]{\ddddot\Pr(#1)}

\title{\Large\vspace*{-42pt}
 The Electrical Response Matrix of a Regular $2n$-gon
} 
\author{Nathaniel D. Blair-Stahn\thanks{University of Washington} \and David B. Wilson\thanks{Microsoft Research}}
\date{}

\begin{document}

\maketitle
\vspace*{-22pt}

\begin{abstract}
  Consider a unit-resistive plate in the shape of a regular polygon
  with $2n$ sides, in which even-numbered sides are wired to
  electrodes and odd-numbered sides are insulated.  The response
  matrix, or Dirichlet-to-Neumann map, allows one to compute the
  currents flowing through the electrodes when they are held at
  specified voltages.  We show that the entries of the response matrix
  of the regular $2n$-gon are given by the differences of cotangents
  of evenly spaced angles, and we describe some connections with the
  limiting distributions of certain random spanning forests. 
\end{abstract}

\section{Introduction} \label{intro sec}

Let $n\ge 2$ be an integer, and let $\poly$ be a regular polygon
with $2n$ sides, centered at the origin in $\C$, with
the midpoint of the $j\th$ side of $\poly$ located
on the unit circle at $e^{i\pi j/n}$ for $1\le j\le 2n$.  We
imagine $\poly$ to be made of a unit resistive material, and we wish
to determine how much current will flow through $\poly$ when
the even-numbered sides of $\poly$ are wired to electrodes at specified
voltages while the odd-numbered sides of $\poly$ are kept insulated.

We will call the even-numbered (wired) sides of $\poly$
\textbf{nodes} and the odd-numbered (insulated) sides \textbf{free
edges}. Current may pass through the nodes, but no current is
allowed to enter or exit $\poly$ through the free edges.  We number
the $n$ nodes so that node $j$ corresponds to side $2j$ of $\poly$.
For each $j$, the $j\th$ node of $\poly$ is wired to an electrode
and held at some specified voltage $v_j$. Each \textbf{voltage
configuration} $\vec v = (v_1,\ldots,v_n)$ on the nodes results in
some \textbf{current output} $\vec I = (I_1,\ldots,I_n)$, where
$I_j$ is the current \emph{entering} $\poly$ through node $j$.

This problem can be rephrased in terms of the solution to a mixed
Dirichlet-Neumann boundary problem on the domain $\poly$:  Given
constants $v_1,\ldots, v_n$, there is a unique continuous function
$V$ (the electric potential) on $\poly$ which is harmonic on the
interior, equals $v_j$ on side $2j$, and has zero normal derivative
on the free edges (see e.g., \cite[pg.~452 or pg.~445]{\GM}).  The
current $I_j$ entering $\poly$ through node $j$ is defined to be the
flux of the electric field $E_V = -\nabla V$ through side $2j$ into
$\poly$.

Since the current output $\vec I$ depends linearly on the voltage configuration
$\vec v$, there is an \mbox{$n\times n$} matrix $\Lambda$ mapping the space of
voltage configurations to the space of current outputs. $\Lambda$ is
called the \textbf{response matrix} or \textbf{Dirichlet-to-Neumann
map} for the ``resistor network'' represented by $\poly$.  The
Dirichlet-to-Neumann map $\Lambda$ is defined analogously for more general
resistive domains or for resistor networks represented by finite
graphs (see e.g.\ \cite{\CdV}, \cite{\IM}, \cite{\CIM}
for background).
Since, for any resistor network, the $k\th$ column of $\Lambda$ is
the action of $\Lambda$ on the $k\th$ standard basis vector,
$\Lambda_{j,k}$ will be the current entering the network ($\poly$ in
our case) through node $j$ when node $k$ is held at 1 volt and all
other nodes are held at 0 volts.  Observe that $\Lambda_{j,k}$ will
be positive if $j=k$ and negative otherwise.  Although it is not
obvious from this description, the response matrix for a general
electrical network is symmetric (see \cite{\CdV}).
In this
paper we compute the response matrix $\Lambda$ for~$\poly$ with
the boundary conditions described above:

\begin{thm} \label{main thm}
The response matrix $\Lambda$ for the regular $2n$-gon with
alternating wired/free boundary conditions has entries given by
$$
\Lambda_{j,k} = \frac{\cot \left[\frac{\pi}{n} \left(j-k+\frac{1}{2}\right)\right]
        - \cot \left[\frac{\pi}{n} \left(j-k-\frac{1}{2}\right)\right]}{n}.
$$
\end{thm}

The response matrix is closely connected to the distribution of random
``groves'' (a generalization of spanning trees) in a resistor network
\cite{\KW}.  In \S~\ref{groves/polynomials sec} we give some background
on this connection with groves, which is what initially led us to study the response matrix
$\Lambda$, and we briefly compare the grove model with an analogous model based on percolation.
We then discuss a purely algebraic approach (based on groves) that Kenyon and Wilson
\cite[\S~5.2]{\KW} used to compute the response matrix
for the regular $2n$-gon in the cases $n=3$ and $n=4$.  This algebraic
approach, however, is not easily adapted to general $n$. We prove Theorem~\ref{main thm}
in \S~\ref{diag sec}, using a combination
of algebraic and analytic methods.

\begin{figure}[bhp]
\psfrag{->}[cc][Bc][1][0]{$\to$}
\centerline{\includegraphics[width=.84\textwidth]{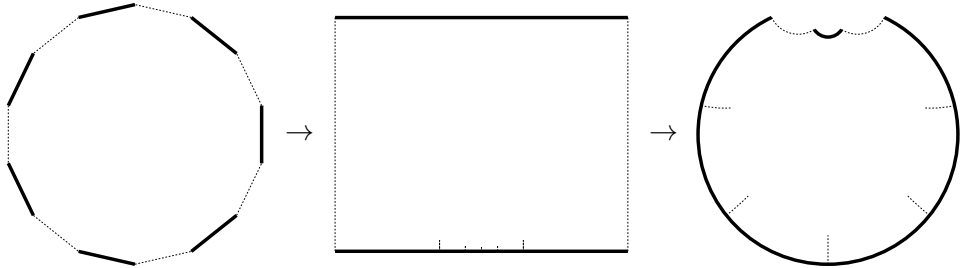}}
\caption{On the left is the regular $2n$-gon ($n=7$) with alternate sides
  wired and free (insulated).  The regular $2n$-gon is conformally
  mapped to a rectangle so that one wired side goes to the top and the
  remaining wired sides go the bottom of the rectangle, while two free
  sides get mapped to the sides of the rectangle and the remaining
  free sides get mapped to vertical slits.  Entries of the response
  matrix are given by the lengths of the wired sides in the rectangle
  divided by the height of the rectangle.  A M\"obius transformation
  maps the slit rectangle to a subset of the unit disk.  Each wired
  side and each free side is mapped to an arc of a circle that passes
  through the top of the disk, and $n-1$ of the wired sides each get
  mapped to arcs of the unit circle of length $2\pi/n$.  }
\label{hept}
\end{figure}

We mention that it is known how to compute the response matrix by
using the Schwarz-Christoffel formula to map the polygon to a
rectangle with vertical slits that correspond to the free edges,
as shown in Figure~\ref{hept}.
A general $2n$-gon may be conformally mapped to a rectangle so that
one wired side (say side $j$) gets mapped to the top of the
rectangle, the adjacent free sides get mapped to the sides of the
rectangle, the remaining wired sides get mapped to intervals of the
bottom side of the rectangle, and the remaining free sides get
mapped to vertical slits.  (\cite[\S~5.2]{\KW} includes a discussion
of these maps.) The current response $\Lambda_{j,k}$ is just the
ratio of the length of the image of the $k\th$ side to the height of
the rectangle. Without going further into details, we mention that
in the example of the regular octagon ($n=4$), this approach yields
$$
\Lambda_{j,j+2} = - \frac{\displaystyle\int_{x_4}^{x_5}\frac{w^2-b^2}{\prod_{\ell=1}^8 (w-x_\ell)^{1/2}}\,dw}{\displaystyle\int_{x_7}^{x_8}\frac{w^2-b^2}{\prod_{\ell=1}^8 (w-x_\ell)^{1/2}}\,dw} \ \ \ \ \text{where} \ \ \ \
b^2 = \frac{\displaystyle\int_{x_3}^{x_4}\frac{w^2}{\prod_{\ell=1}^8 (w-x_\ell)^{1/2}}\,dw}{\displaystyle\int_{x_3}^{x_4}\frac{1}{\prod_{\ell=1}^8 (w-x_\ell)^{1/2}}\,dw}$$
and $x_\ell = \cot((1/2-\ell)\pi/8)$.
It is not at first obvious that this should simplify to
$$\Lambda_{j,j+2}=1/2-1/\sqrt2.$$

\remark
Since the current responses for the regular polygon 
are differences of cotangents, the horizontal positions of
the vertical slits and the sides of the rectangle are given by
cotangents of evenly spaced angles.  If we view the slit rectangle as
being embedded in the upper half plane with the rectangle's bottom
edge on the real axis, then there is a M\"obius transformation of
the upper half plane to the unit disk so that images of the bases of
the vertical slits (and the two bottom corners) are evenly spaced on
the unit circle (Figure~\ref{hept}).

\section{The response matrix and random groves} \label{groves/polynomials sec}

\begin{figure}[b!]
\psfrag{rX1 rectangle}[cc][Bc][1][90]{$r\times 1$ rectangle}
\psfrag{regular hexagon}[cc][Bc][1][90]{regular hexagon}
\psfrag{Cardy, Ziff, Smirnov}[tc][tc][1][0]{\parbox{2.6in}{\center $\Pr[\sigma\!=\!12]\xrightarrow[\varepsilon\to 0]{}  \frac{\Gamma(\frac23)\, s^{1/3}}{\Gamma(\frac43)\Gamma(\frac13)}\, _2F_1(\frac23,\,\frac13;\,\frac43;\,s)$ \\ where $s = \left(\frac{\sum_{n\in\Z} e^{-\pi r (n+1/2)^2}}{\sum_{n\in\Z} e^{-\pi r n^2}}\right)^4$ \\ \cite{MR92m:82048,MR1851632,ziff2} }}
\psfrag{Kirchhoff}[tc][tc][1][0]{$\Pr[\sigma=12]=\frac{-\Lambda_{1,2}}{1-\Lambda_{1,2}}\xrightarrow[\varepsilon\to 0]{}\frac{1}{1+r}$\ \ \cite{Kirchhoff}}
\psfrag{Dubedat}[tc][tc][1][0]{\parbox{3.3in}{\vspace*{-8pt}\center $\Pr[\sigma=13|2]\xrightarrow[\varepsilon\to 0]{}\frac{3^{3/2} \Gamma(\frac23)^9}{2^{7/3} \pi^5} \, _3F_2(1,\frac56,\frac56;\, \frac32,\frac32;\, 1)$ \\ \ \ \ \ \ \ \ \cite{MR2253875} }}
\psfrag{KW}[tc][tc][1][0]{$\begin{aligned}\Pr[\sigma=13|2]&=\textstyle\frac{-\Lambda_{1,3}}{\genfrac{}{}{0pt}{}{1-\Lambda_{1,2}-\Lambda_{1,3}-\Lambda_{2,3}+}{\Lambda_{1,2}\Lambda_{1,3}+\Lambda_{1,2}\Lambda_{2,3}+\Lambda_{1,3}\Lambda_{2,3}}}\\&\xrightarrow[\varepsilon\to 0]{} 2/\sqrt3 -1\ \ \cite{\KW}\end{aligned}$}
\centerline{\includegraphics[width=.85\textwidth]{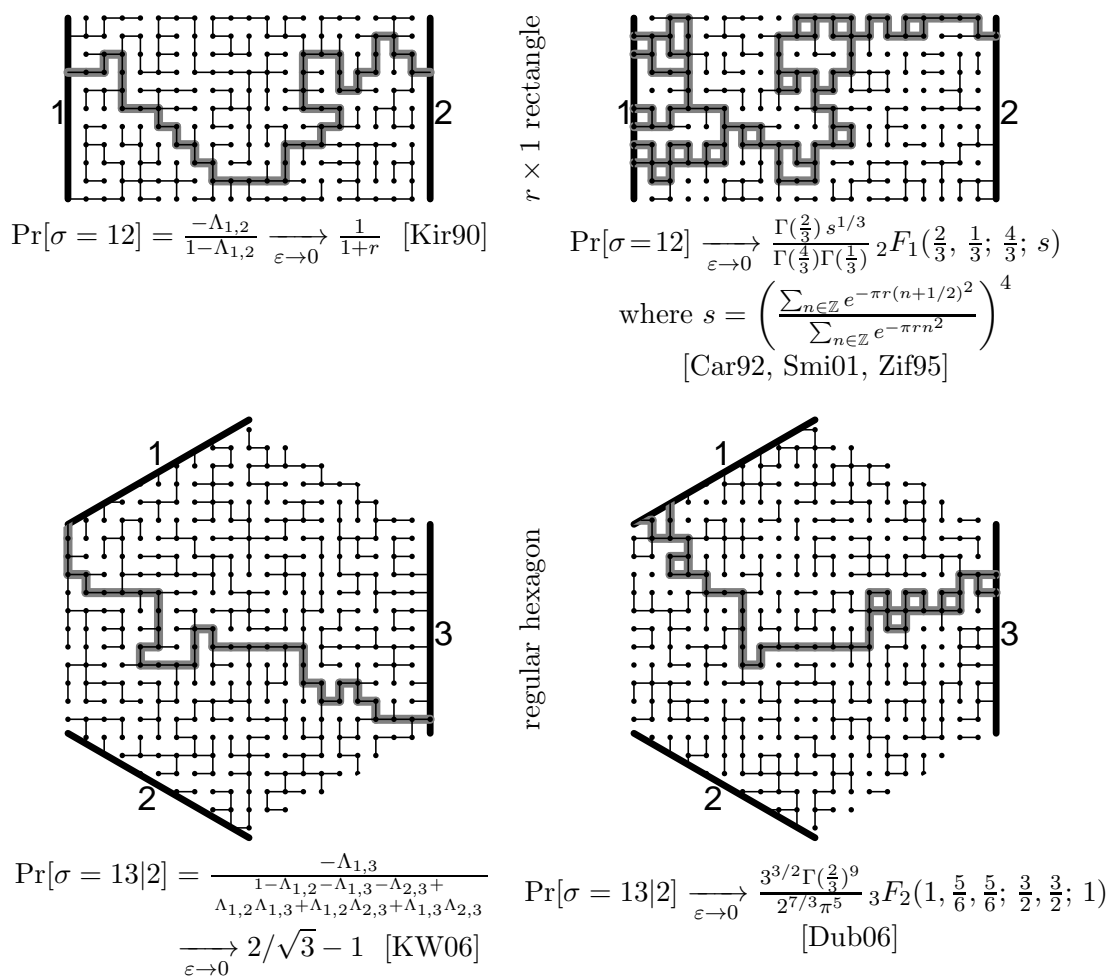}}
\centerline{}
\centerline{}
\caption{
  Examples of grove (left) and percolation (right) configurations on
  regions with two (upper) or three (lower) nodes.  Random groves and
  percolation are both important special cases of the
  Fortuin-Kasteleyn random cluster model from statistical physics
  \cite{FK}.  The paths connecting the nodes are highlighted, and when
  the grid spacing~$\varepsilon$ tends to $0$, the
  probability that the nodes are connected in a given way is indicated
  below the panels.  (The percolation formulas are rigorous for site
  percolation on the triangular lattice, and are conjectured to hold
  for the bond percolation shown here.)}
\label{perc-grove}
\end{figure}

We give here some background on the relation between the response
matrix of a graph and random ``groves'' of that graph.  A grove is a
forest, i.e., an acyclic collection of edges of the graph, such that
every constituent tree of the forest contains at least one of a
special set of distinguished vertices, which are called nodes.  The
upper-left panel of Figure~\ref{perc-grove} shows a grove on a graph
with two nodes (labeled $1$ and $2$), and the lower-left panel shows
a grove on a graph with three nodes (labeled $1$, $2$, and $3$).  The
first grove consists of one tree, which contains the two nodes $1$ and
$2$, and the path connecting the nodes is highlighted.  The second grove
consists of two trees, one of which contains nodes~$1$ and~$3$, while
the other contains just node~$2$.
\enlargethispage{12pt}
\newpage

Groves of ``circular planar graphs'' arise naturally in combinatorics
\cite{MR2097339} and statistical physics \cite{KW-polynomial}.  A
graph with distinguished nodes is ``circular planar'' if it embeds in
the plane and each of the nodes lies on the outer face; the nodes are
numbered in a counterclockwise order.  If a uniformly random grove is
chosen, this defines a random (noncrossing) partition $\sigma$ on the
set of nodes, and we are interested in the probability distribution of
this random partition.

Another natural probability distribution on node partitions arises
from percolation on the graph (where each edge is included
independently with probability $p$); this is illustrated in the
right-hand panels of Figure~\ref{perc-grove} for $p=1/2$.
There has been significant recent work in studying crossing events for
percolation such as those shown in Figure~\ref{perc-grove}, and it is
interesting to compare the crossing event probabilities in these two
models.

For graphs with two nodes, the partition
probabilities for groves follow from Kirchhoff's formula \cite{Kirchhoff}:
\begin{equation}\label{Kirchhoff formula}
\frac{\text{\# spanning trees}}{\text{\parbox{1.7in}{\small\ \ \#
      two-tree forests with \\[-2pt] nodes $1$ and $2$ disconnected}}}
= \frac{\Pr[\text{grove partition is $12$}]}{\Pr[\text{grove partition
    is $1|2$}]} = \frac{1}{R_{1,2}} = -\Lambda_{1,2},
\end{equation}
where
$R_{1,2}$ is the electrical resistance between nodes $1$ and $2$ when
the graph is viewed as a resistor network in which each edge has
unit resistance.  (The relation
$\Lambda_{1,2}=-1/R_{1,2}$ holds for two nodes only; when there are
more nodes, the response matrix entries are more complicated functions
of the pairwise resistances.)

In the case where the graph is a fine grid restricted to an $r \times
1$ rectangular region, with the two nodes being extra vertices
corresponding to the left and right edges of the rectangle (as
illustrated in Figure~\ref{perc-grove}), then $R_{1,2}=r$, so the
probability that the grove partition $\sigma$ is $12$ is simply
$\Pr[\sigma=12]=1/(1+r)$.  If the fine grid is restricted to a region
different from a rectangle, with two nodes each occupying some
fraction of the boundary, then the region may be conformally mapped to
a rectangle of some aspect ratio $r$, with the two nodes getting
mapped to the left and right edges.  When the grid becomes very fine,
the resistance between the nodes in the orginal domain converges to
$r$, so in this limit $\Pr[\sigma=12]\to 1/(1+r)$.

The corresponding formula for critical percolation (where each edge
occurs independently with probability $1/2$) was deduced by Cardy
\cite{MR92m:82048} (see upper-right panel of Figure~\ref{perc-grove})
using exact but nonrigorous methods.  (Ziff adapted Cardy's formula,
which was given for the upper half plane, to rectangular regions
\cite{ziff2}, and Smirnov gave a rigorous proof for the related model
of site percolation \cite{MR1851632}.)

More recently these boundary crossing events have been studied in
regions with more nodes \cite{MR2253875,\KW}; this is illustrated in
the lower-right (for percolation) and lower-left (for groves) panels
of Figure~\ref{perc-grove}.  In the case of groves, Kenyon and Wilson
\cite{\KW} show how to compute the boundary partition probabilities
for any number of nodes in terms of the entries of the response matrix
$\Lambda$ of the graph when viewed as an electrical network.
It is convenient to abbreviate
$$
\pu{\sigma} =  \frac{\Pr[\text{grove partition is $\sigma$}]}{\Pr[\text{grove partition is $1|2|\dots|n$}]}.$$
For graphs with three nodes the analogues of Kirchhoff's formula \eqref{Kirchhoff formula} are
\begin{equation}\label{3 node formulas}
\begin{matrix}
\pu{123}=\Lambda_{1,2}\Lambda_{1,3}+\Lambda_{1,2}\Lambda_{2,3}+\Lambda_{1,3}\Lambda_{2,3},\\
\pu{1|23} = -\Lambda_{2,3},\ \ \ \ \ \ 
\pu{2|13} = -\Lambda_{1,3},\ \ \ \ \ \ 
\pu{3|12} = -\Lambda_{1,2},\ \ \ \ \ \ 
\pu{1|2|3}= 1.
\end{matrix}
\end{equation}
More generally, Kenyon and Wilson \cite{\KW} proved that for any
circular planar graph with any number $n$ of nodes on the outer face,
if $\sigma$ is any noncrossing partition on $\{1,\dots,n\}$, then
$\pu{\sigma}$ is a polynomial in the entries of the response matrix
$\Lambda$ and can be computed explicitly.

It is interesting to see what these general formulas give for some
nice special cases.  For example, Dub\'edat \cite{\Dub} computed the partition
probabilities for (site) percolation when the region is the regular
hexagon with 3 nodes along alternate sides of the hexagon.  (The
lower-right panel of Figure~\ref{perc-grove} illustrates this for bond
percolation.)  To carry out a similar computation for groves, we need
to calculate the response matrix for the regular hexagon with
alternating free and wired boundary conditions on its faces.

In the introduction we mentioned an algebraic approach for computing
the response matrix $\Lambda_{j,k}$ for a regular $2n$-gon.  This
approach \cite[\S~5.2]{\KW} is based on enumerating trees in a random
grove; we briefly describe this tree-enumeration approach together
with some associated open problems.

We consider random groves on a very fine grid restricted to the
regular $2n$-gon, with $n$ nodes ``wired'' to every other side of the
polygon (as in Figure~\ref{perc-grove}).
Using the formulas from \cite{\KW}, we may express the polynomial
$$
P_n(q) = \sum_{t=1}^n \frac{\Pr[\text{$t$ trees in grove}]}{\Pr[\text{$n$ trees in grove}]}\, q^{t-1}
$$
in terms of the response matrix of the graph, which in the scaling
limit approaches the response matrix $\Lambda$ of the regular
$2n$-gon. Since (by symmetry considerations) the response matrix for
the regular $2n$-gon is circulant, we may define
$\Lambda_{|j-k|}=\Lambda_{j,k}$ (where indices are identified mod~$n$).
From the above formulas \eqref{3 node formulas} we see
that, in the limit, $P_3(q) = 1 - 3\Lambda_1 q + 3 \Lambda_1^2 q^2$.
For the regular octagon it turns out that
 $$P_4(q) = 1-q(4 \Lambda_1+2\Lambda_2)+q^2(6\Lambda_1^2+8\Lambda_1 \Lambda_2+2\Lambda_2^2) -q^3(4\Lambda_1^3+8\Lambda_1^2\Lambda_2+4\Lambda_1\Lambda_2^2).$$
 
 For each grove of a circular planar graph with $n$ nodes there is a
 dual grove on the dual graph, which contains the duals of edges not
 contained in the primal grove.  The number of trees in a grove plus
 the number of trees in its dual grove is $n+1$.
 Since the dual graph of a fine grid restricted to the regular $2n$-gon 
 with alternate wired/free boundary conditions is again a fine grid
 restricted the the regular $2n$-gon with alternate free/wired boundary
 conditions, in the
 limit where the grid is very fine, the probability of seeing $t$
 trees in a random grove equals the probability of seeing $n-t+1$
 trees.  Hence the coefficients of the polynomial $P_n(q)$ form a
 palindrome.  For $n=3$ this implies that $3\Lambda_1^2=1$, or
 $\Lambda_1=-1/\sqrt3$, determining the response matrix.  For $n=4$
 there are several pairs $(\Lambda_1,\Lambda_2)$ which make $P_4(q)$ a
 palindrome, but only one in which $\Lambda_1$ and $\Lambda_2$ are
 both negative ($\Lambda_1=-1/2$, $\Lambda_2=1/2-1/\sqrt 2$), which is
 enough to determine the response matrix.  For general $n$ the fact
 that $P_n(q)$ is a palindromic polynomial generates enough
 constraints to limit the coefficients of $\Lambda$ to a
 zero-dimensional algebraic variety, but it is not clear that there
 will always be a unique negative solution for the $\Lambda$'s, nor is
 it clear how to obtain the solution for general $n$ using this
 algebraic approach.
 
 However, in the next section we explicitly compute the response
 matrix for the polygon~$\poly$ by other means, proving
 Theorem~\ref{main thm}. Combining this result with the formulas from
 \cite{\KW}, we may write down the first few polynomials $P_n(q)$:
\begin{align*}
P_2(q) &= 1 + q \\
P_3(q) &= 1 + \sqrt{3} q + q^2\\
P_4(q) &= 1 + (1+\sqrt2) q + (1+\sqrt2) q^2 + q^3\\
P_5(q) &= 1 + \sqrt{5+2\sqrt5} q + (2+\sqrt{5}) q^2 + \sqrt{5+2\sqrt5} q^3 + q^4\\
P_6(q) &= 1 + (2+\sqrt3) q + (3+2\sqrt3) q^2 + (3+2\sqrt3) q^3 + (2+\sqrt3) q^4 + q^5
\end{align*}
The constant term is of course always $1$.  Referring to \cite{\KW},
the linear term is $-\sum_{j<k} \Lambda_{j,k}$, which simplifies to
$\cot(\pi/(2n))$.  It would be interesting to better understand the
polynomials $P_n$, such as for example the approximate distribution of
the number of trees for large~$n$, or the asymptotics of $P_n(1) =
1/\Pr[\text{random grove is a single tree}]$.

\section{Computing $\Lambda$} \label{diag sec}

This section is devoted to proving Theorem~\ref{main thm}.  In
\S~\ref{evals/evecs sec} we identify the
eigenvectors and eigenvalues of $\Lambda$ using symmetry
considerations, thereby finding a diagonalization of~$\Lambda$.  Then
in \S~\ref{matrix mult sec} we compute $\Lambda$ by performing a
matrix multiplication to change from the eigenbasis back to the
standard basis.

\subsection{Complex potential and the diagonalization of $\Lambda$} \label{evals/evecs sec}

The key to finding the eigenvectors and eigenvalues of $\Lambda$ is to
introduce complex electric potentials and currents in order to exploit
the symmetry of the polygon $\poly$.  A complex-valued potential $V$
on $\poly$ can be thought of as two separate real potentials, one from
the real part of $V$ and the other from the imaginary part.  The
electric field $E_V$ now takes values in $\C^2$ rather than $\R^2$,
and can be thought of as carrying separate real and imaginary
currents.  The current output $\vec I$ will now be a \emph{complex}-linear
function of the voltage configuration $\vec v$. Thus we can view $\Lambda$
as a complex-linear transformation of $\C^n$ whose restriction to
$\R^n$ yields the expected real current outputs.

Let $\omega = e^{2\pi i/n}$, and for $1\le k\le
n$, define the voltage configuration $\vk$ by $(\vk)_j = \omega^{jk}$.

\begin{lem} \label{evals/evecs lem}
Let $1\le k\le n$. Then $\vk = (\omega^k,\omega^{2k},\ldots,\omega^{(n-1)k}, 1)$ is an eigenvector of $\Lambda$, and the corresponding eigenvalue $\lambda_k$ is the current entering $\poly$ through node $n$ under the voltage configuration $\vk$.
\end{lem}

\begin{proof}
Let $\sigma: \C^n\rightarrow \C^n$ denote the function which
cyclically shifts the components of a vector to the left: $\sigma
(v_1,v_2,\ldots,v_{n-1},v_n) = (v_2,v_3,\ldots, v_n,v_1)$. If $\vec v$ is
a voltage configuration on $\poly$, $\sigma \vec v$ is the voltage
configuration obtained by replacing the voltage on node~$j$ with the
voltage on node $j+1$ (where the indices are identified mod~$n$).
Because $\poly$ is rotationally symmetric, the resulting currents
will likewise be rotated clockwise by one node.  That is,
\begin{equation*}
\Lambda (\sigma \vec v) = \sigma (\Lambda \vec v) {\rm \ \ for\ all\ } \vec v\in \C^n.
\end{equation*}
On the other hand, for each vector $\vk = (\omega^k,\omega^{2k},\ldots,\omega^{(n-1)k}, 1)$, we have $\sigma \vk = \omega^k \vk$.
If $\Ik=\Lambda \vk$ is the current output resulting from $\vk$, then
we have
\begin{equation*}
\sigma \Ik = \sigma(\Lambda\vk) = \Lambda (\sigma\vk) = \Lambda(\omega^k \vk) = \omega^k \Ik.
\end{equation*}
This implies $(\Ik)_j = \omega^{jk}(\Ik)_n$, so $\Lambda \vk = \Ik =
(\Ik)_n \vk$.
\end{proof}

Next we compute the eigencurrent $\lambda_k = (\Ik)_n$ by
considering the harmonic conjugate of the potential induced by
$\vk$. Since harmonic conjugation defines a real-linear operator --
call it $\harmconj$ -- on real-valued harmonic functions on $\poly$
(modulo constant functions), there is a unique way to extend
$\harmconj$ to a complex-linear operator on complex-valued harmonic
functions on $\poly$. Recall that any harmonic conjugate satisfies
the Cauchy-Riemann equations,
$$
(\harmconj V)_x = -V_y \quad \text{and} \quad (\harmconj V)_y = V_x,
$$
so that $\nabla \harmconj V$ is orthogonal to $\nabla V$ at every
point. The Cauchy-Riemann equations imply the following well-known
result (see e.g., \cite[Section~III.6]{\Gam}),
which we shall use to compute the current $\lambda_k$:
\begin{prop} \label{conj/current prop}
Let $\Omega\subseteq \C$ be simply connected, let $V:\Omega \rightarrow \C$
be harmonic, and let \mbox{$\gamma:[0,1]\rightarrow \Omega$} be a $C^1$
path with $\gamma(0)=a$ and $\gamma(1)=b$.  If $\harmconj V$ is any
harmonic conjugate of $V$, then the current due to $V$ flowing from
left to right across $\gamma$ is equal to $\harmconj V(a) -
\harmconj V(b)$ (where `left' and `right' are defined relative to
$\gamma$'s orientation).
\end{prop}

Now let $\harmconj_0 V$ denote the harmonic conjugate of $V$ satisfying
$(\harmconj_0 V)(0) = i V(0)$, so that $\harmconj_0^2 = -\identity$, and
let $\Vk :\poly\rightarrow \C$ be the potential function induced by the
voltage configuration $\vk$.

\begin{figure}[phtb]
\psfrag{z2}[Bc][Bc][1][0]{$\zeta^{2k}$}
\psfrag{z4}[cr][Br][1][0]{$\zeta^{4k}$}
\psfrag{z6}[Bc][Bc][1][0]{$\zeta^{6k}$}
\psfrag{z8}[cl][Bl][1][0]{$\zeta^{8k}$}
\psfrag{z1}[Bl][Bl][1][0]{$-i\zeta^{k}$}
\psfrag{z3}[Br][Br][1][0]{$-i\zeta^{3k}$}
\psfrag{z5}[Br][Br][1][0]{$-i\zeta^{5k}$}
\psfrag{z7}[Bl][Bl][1][0]{$-i\zeta^{7k}$}
\psfrag{V}[cc][Bc][1][0]{$\Vk$}
\psfrag{HV}[cc][Bc][1][0]{$\harmconj_0 \Vk$}
\centerline{\includegraphics[width=0.5\textwidth]{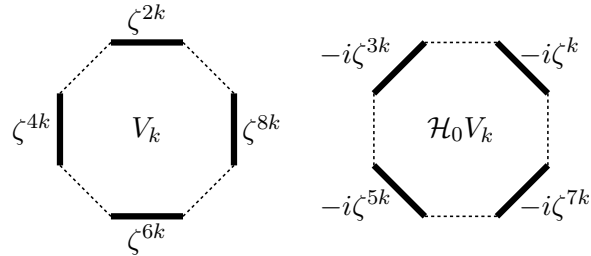}}
\caption{Boundary conditions of $\Vk$ and $\harmconj_0 \Vk$ for
$n=4$, illustrating Lemma~\ref{conj bndry lem}.
Thick solid lines represent wired sides of $\poly$, and thin dashed lines
represent free sides.}
\label{bndry cond fig}
\end{figure}

\begin{lem} \label{conj bndry lem}
  Let $\zeta= e^{i\pi/n}$ (so $\zeta^2=\omega$).  Then $\harmconj_0
  \Vk$ takes the value $-i \zeta^{(2j-1)k}$ on side $2j-1$, for $1\le
  j\le n$, and has zero normal derivative on the even-numbered sides
  of $\poly$ (see Figure~\ref{bndry cond fig}).  Under the
  potential~$\Vk$, the current $\lambda_k$
  entering $\poly$ through side $2n$ is $2\sin (\pi k/n)$.
\end{lem}

\begin{proof}
The fact that the boundary conditions of $\harmconj_0 \Vk$ are
wired where $\Vk$ is free and vice versa is immediate from the
Cauchy-Riemann equations:
Since $\nabla \Vk$ is oriented parallel to the boundary on the free
edges and orthogonal to the boundary on the nodes, the reverse is true for $\nabla
\harmconj \Vk$.  (Boundary issues can be dealt with by using Schwarz
reflection to enlarge the domain, moving the boundary sides to the
interior.)

Recall that rotating the
voltages in the configuration $\vk$ clockwise by one node is the
same as multiplying them by $\omega^k$.  Thus, if $z$ is on an
even-numbered side of $\poly$, then
\begin{equation} \label{rotate potential eqn}
\Vk(\omega z) = \omega^k \Vk(z).
\end{equation}
Now for each $\alpha\in\C$ with $|\alpha|=1$, we define the rotation
operator $R_\alpha(z) = \alpha z$ for $z\in\C$.  With this notation,
\eqref{rotate potential eqn} says that the functions $\Vk\circ
R_\omega$ and $\omega^k \Vk$ agree on the nodes of $\poly$.  Since
they are also both continuous on $\poly$, harmonic on the interior,
and have zero normal derivative on the free edges, they must be
equal, so \eqref{rotate potential eqn} in fact holds for all
$z\in\poly$.

Since $\harmconj$ is linear, \eqref{rotate potential eqn} implies
$ \harmconj_0 \Vk(\omega z) = \omega^k \harmconj_0 \Vk(z) + C $
for some constant $C$.  Setting $z=0$ shows that
$C=(1-\omega^k)\harmconj_0 \Vk(0)=0$ (because $\harmconj_0 \Vk(0)=
i\Vk(0)=0$ if $k\ne n$, and $1-\omega^n=0$).  Thus we have
\begin{equation} \label{rotate conj eqn}
  \harmconj_0 \Vk(\omega z) = \omega^k \harmconj_0 \Vk(z)
\end{equation}
for all $z\in\poly$, so $\harmconj_0 \Vk$ has the same type of
rotational symmetry as $\Vk$.

Let $a_k$ denote the (constant) value of $\harmconj_0 \Vk$ on
side 1.  Combining \eqref{rotate conj eqn} with the fact that
$\harmconj_0 \Vk$ has alternating free/wired boundary conditions
shows that the boundary conditions of the harmonic functions $a_k \Vk$ and
$(\harmconj_0 \Vk) \circ R_\zeta$ agree, so these functions must be equal. Therefore,
\begin{align*} \label{HV/V related eqn}
\harmconj_0 \Vk &= a_k \Vk \circ R_{\zeta^{-1}} \\
\harmconj_0^2 \Vk &= a_k \harmconj_0(\Vk \circ R_{\zeta^{-1}}) \\
-\Vk &= a_k (\harmconj_0 \Vk) \circ R_{\zeta^{-1}} \\
-\Vk &= a_k^2 \Vk \circ R_{\zeta^{-2}} \\
-\Vk &= a_k^2 \omega^{-k} \Vk,
\end{align*}
which shows that $a_k = \pm i\zeta^k$ (note that \eqref{rotate
potential eqn} was used in the last step).

If $k=n$ then $\Vk \equiv 1$, which implies
$\harmconj_0 \Vk \equiv i$, so $a_n = i = -i\zeta^n$ and $\lambda_n = 0$.
If $1\le k < n$, we use Proposition~\ref{conj/current prop} to compute the current
$\lambda_k$, up to a choice of sign:
$$
\lambda_k = \harmconj_0 \Vk (e^{i\pi/2n}) - \harmconj_0 \Vk
(e^{-i\pi/2n}) = a_k-\omega^{-k}a_k = \pm i\zeta^k \mp i\zeta^{-k} =
\mp 2\sin(\pi k/n).
$$
To determine the sign, note that
since $\Vk$ is harmonic, we have $\Vk(z)={\mathbb E}[\Vk(z_T)]$ for
all $z\in\poly$, where $z_t$ is a standard Brownian motion started
at $z_0=z$ which is reflected off the odd sides and absorbed at the
even sides, and $T$ is the absorption time.  On the even-numbered
sides, $\Re \Vk(z)$ is maximized (with value 1) on side $2n$.  Thus
$\Re \Vk(z)<1$ for $z$ in the interior of the polygon $\poly$, so
the real part of the electric field on side $2n$ points into~$\poly$.
Hence the real part of the current entering $\poly$ on side
$2n$ must be nonnegative, so it is $+2\sin(\pi k/n)$, and $a_k =
-i\zeta^k$.
\end{proof}

\subsection{Recovering $\Lambda$ from its diagonalization} \label{matrix mult sec}

Putting together Lemmas~\ref{evals/evecs lem} and
\ref{conj bndry lem} we have
\begin{thm} \label{diag thm}
The response matrix $\Lambda$ satisfies $\Lambda W = WD$, where $W$ is the matrix of eigenvectors given by $W_{j,k}=e^{2\pi i jk/n}$, and $D$ is the diagonal matrix of eigenvalues with $D_{k,k} = \lambda_k = 2\sin(\pi k/n)$.
\end{thm}


We use the following lemma to compute $\Lambda$ from Theorem~\ref{diag thm}.  Recall that $\zeta = e^{i\pi/n}$.

\begin{lem} \label{symmetric sums lem}
If $m$ is any integer, then
$$
\sum_{\ell=1}^n \zeta^{(2m\pm 1)\ell}
    = -1 + i\, \textstyle{\cot\left[\frac{\pi}{n}\left(m\pm \frac{1}{2}\right)\right]}.
$$
\end{lem}

\begin{figure}[phtb]
    \psfrag{Beta}[Br][Br][1][0]{$\beta$}
    \psfrag{Beta'}[cl][Bl][1][0]{$\ \beta' = e^{-i\theta}\beta$}
    \psfrag{t}[cc][Bc][1][0]{$\,\theta$}
    \psfrag{cott }[bc][Br][1][90]{$\overbrace{\hspace{111.962bp}}^{\cot\theta}$}
    \psfrag{z}[cl][Bl][1][0]{$\,\zeta=e^{i\pi/n}$}
    \psfrag{-1}[tr][tr][1][0]{$-1$}
\centerline{\includegraphics{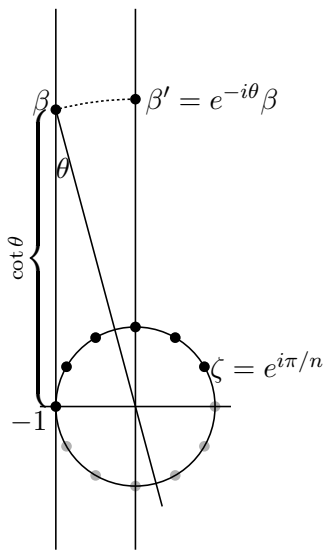}}
\caption{Illustration of the proof of Lemma~\ref{symmetric sums lem}
with $n=6$ and $m=0$. The first five powers of $\zeta$ are symmetric
about the imaginary axis, so the sum $\beta = \sum_{\ell=1}^6
\zeta^\ell$ lies on the line $\Re z = -1$. The first six powers of $\zeta$ are symmetric about the line $\arg z = \theta +\pi/2$, where $\theta= \pi/12$, so $\beta$ also lies on this line (equivalently, $\beta'= e^{-i\theta} \beta$ is imaginary).}
\label{cot fig}
\end{figure}

\begin{proof}
Let $b_\ell=\zeta^{(2m\pm 1)\ell}$  and let $\beta = \sum_{\ell=1}^n
b_\ell$.  We first show that $\Re \beta = -1$. Since $b_{n-\ell} =
-\overline{b_\ell}$, the real parts of $b_\ell$ and $b_{n-\ell}$
add to 0
 (or the real part equals 0 if $\ell=n/2$).  Pairing up the terms in
 this way, we see $\Re \sum_{\ell=1}^{n-1} b_\ell = 0$ and hence $\Re\beta = \Re b_n=-1$.

Now let $\theta =
\frac{\pi}{n} \left(m\pm \frac{1}{2}\right)$, and define $\beta' =
e^{-i\theta}\beta$.  We will use the same trick as above to show
that $\beta'$ is imaginary, then show that this implies $\Im
\beta = \cot \theta$ (see Figure \ref{cot fig}). Note that $\beta' =
\sum_{\ell=1}^n c_\ell$, where $c_\ell = e^{-i\theta}b_\ell =
e^{\frac{i\pi}{2n} (2m\pm 1) (2\ell-1)}$. Now,
$$
c_{n+1-\ell} = e^{\frac{i\pi}{2n} (2m\pm 1) (2n-2\ell+1)}
    = e^{\pm i\pi}e^{-\frac{i\pi}{2n} (2m\pm 1) (2\ell-1)}
    = -\overline{c_\ell},
$$
so $c_\ell$ and $c_{n+1-\ell}$ are symmetric about the imaginary
axis. It follows that $\beta'$ is imaginary and hence $\arg \beta =
\theta + \arg \beta' = \theta \pm \frac{\pi}{2}$.  Finally, observe
that
\begin{align*}
\Im \beta = \Re \beta \cdot \tan (\arg\beta)
    = -1\cdot \textstyle{\tan\left(\theta \pm \frac{\pi}{2}\right)}
    &= \cot \theta
    = \textstyle{\cot\left[\frac{\pi}{n}\left(m\pm \frac{1}{2}\right)\right].} \qedhere
\end{align*}
\end{proof}

\extraline

\begin{proof}[Proof of Theorem~\ref{main thm}]
First we observe that the matrix $W$ in Theorem~\ref{diag thm} is
invertible with $W^{-1} = \frac{1}{n}W^*$. (One can easily verify
this directly, or simply notice that $\frac{1}{\sqrt n}W$ is the
inverse discrete Fourier transform matrix, which is unitary.) We
need to show that the entries of the matrix $WDW^{-1}$ agree with
the formula for $\Lambda$ given in Theorem~\ref{main thm}.   It
follows from Theorem~\ref{diag thm} that $(WD)_{j,k} = W_{j,k}
D_{k,k} = \zeta^{2jk}\lambda_k = \zeta^{2jk}(\zeta^k-\zeta^{-k})/i$.
Since $(W^{-1})_{j,k}=\frac{1}{n}\overline{W_{k,j}}=\zeta^{-2jk}/n$,
we have
\begin{align*}
\Lambda_{j,k}
    =(WDW^{-1})_{j,k} = \sum_{\ell=1}^n (WD)_{j,\ell}(W^{-1})_{\ell,k}
        &= \frac{1}{in} \sum_{\ell=1}^n \zeta^{2j\ell}(\zeta^\ell-\zeta^{-\ell})\zeta^{-2\ell k}\\
        &= \frac{1}{in} \sum_{\ell=1}^n \left[\left(\zeta^{2(j-k)+1}\right)^\ell -
            \left(\zeta^{2(j-k)-1}\right)^\ell\right]\\
        &= \frac{\cot \left[\frac{\pi}{n} \left(j-k+\frac{1}{2}\right)\right]
            - \cot \left[\frac{\pi}{n}
            \left(j-k-\frac{1}{2}\right)\right]}{n},
\end{align*}
where the last equality follows from Lemma~\ref{symmetric sums lem}.
\end{proof}

\extraline

\noindent \textbf{\Large Acknowledgements}

We thank Richard Kenyon for helpful discussions,
in particular for suggesting the use of complex voltages.

\bibliographystyle{halpha}
\bibliography{bc,other}

\end{document}